\documentclass[11pt]{article}

\usepackage{amsfonts}
\usepackage{amssymb}
\usepackage{graphicx}
\usepackage{cite}
\usepackage{color}

\newcommand{\matt}[1]{\left[ \matrix{#1} \right]}
\def\scr#1{{\cal #1}}
\newcommand{\R}{\mathbb{R}}
\def\eq#1{\begin{equation}#1\end{equation}}
\def\rep#1{(\ref{#1})}
\newcommand{\bbb}{\mathbb}
\newtheorem{theorem}{Theorem}
\newtheorem{lemma}{Lemma}
\newtheorem{remark}{Remark}
\newtheorem{proposition}{Proposition}
\newtheorem{corollary}{Corollary}
\newtheorem{assumption}{Assumption}

\def\qed{ \rule{.1in}{.1in}}

\begin{document}

\title{Exponential Convergence of the Discrete-Time Altafini Model\thanks{
This research was supported in part by AFOSR MURI Grant FA 9550-10-1-0573.
J. Liu, X. Chen, T. Ba\c{s}ar, and M.-A. Belabbas are with the Coordinated Science Laboratory,
University of Illinois at Urbana-Champaign (\texttt{\{jiliu, xdchen, basar1, belabbas\}@illinois.edu}).
}
}

\author{Ji Liu  \hspace{.3in}  Xudong Chen  \hspace{.3in}  Tamer Ba\c{s}ar  \hspace{.3in}  Mohamed Ali Belabbas}
\maketitle


\begin{abstract}

This paper considers the discrete-time version of Altafini's model for opinion dynamics
in which the interaction among
a group of agents is described by a time-varying signed digraph.
Prompted by an idea from \cite{lift}, exponential convergence of the system
is studied using a graphical approach.
Necessary and sufficient conditions for exponential convergence with respect to
each possible type of limit states are provided.
Specifically, under the assumption of repeatedly jointly strong connectivity, it is shown that (1)
a certain type of two-clustering will be reached exponentially fast for almost all initial conditions
if, and only if, the sequence of signed digraphs is repeatedly jointly structurally balanced
corresponding to that type of two-clustering;
(2) the system will converge to zero exponentially fast for all initial conditions
if, and only if, the sequence of signed digraphs is repeatedly jointly structurally unbalanced.
An upper bound on the convergence rate is also provided.

\end{abstract}

\section{Introduction}

Over the past decade, there has been considerable interest in developing algorithms
intended to cause a group of multiple agents to reach a consensus in a distributed manner
\cite{Ts3,Ts4,vicsekmodel,reza1,luc,ReBe05,blondell,survey,basar,reachingp1,hendrickx,cut,touri2014,cdc14}.
Consensus processes have a long history in social sciences and are closely related to opinion dynamics \cite{french}.
Probably the most well-known opinion dynamics is the DeGroot model
which is a linear discrete-time consensus process \cite{degroot}.
Recently, quite a few models have been proposed for opinion dynamics,
including the Friedkin-Johnsen model \cite{johnsen1,johnsen,magazine}, the Hegselmann-Krause model \cite{krause1,krause,etesami}, the Deffuant-Weisbuch model \cite{dw},
and the DeGroot-Friedkin model \cite{jia,Bullo5,Bullo4,acc15}.
A particularly interesting opinion dynamics model, which was first proposed by Altafini \cite{plos,altafini}
and can be viewed as a more general linear consensus model,
has received increasing attention lately \cite{ming1,ming2,mingtac,xia,xia1,weiguo,modulus,lift,signcdc13}.

The Altafini model deals with a network of $n>1$ agents
and the constraint that each agent is able to receive information only from its ``neighbors''.
Unlike the existing models for opinion dynamics or consensus, the
neighbor relationships among the agents are described by a time-dependent, {\em signed}, digraph (or directed graph)
in which vertices correspond to agents, arcs (or directed edges) indicate the directions of information flow,
and, in particular, the signs represent the social relationships between neighboring agents
in that a positive sign means friendship (or cooperation) and a negative sign indicates antagonism (or competition).
Each agent $i$ has control over a time-dependent state variable $x_i(t)$ taking values in $\R$,
which denotes its opinion on some issue.
Each agent updates its opinion based on its own current opinion, the current opinions of its current neighbors,
and its relationships (friendship or antagonism) with respect to its current neighbors.
Specifically,
for those neighbors with friendship, the agent will trust their opinions;
for those neighbors with antagonism, the agent will not trust their opinions and, instead,
the agent will take the opposite of their opinions in updating,
which is the key difference between the Altafini model and other opinion dynamics models.

The continuous-time Altafini model has been studied in \cite{signcdc13,ming1,ming2,mingtac},
and papers \cite{xia,xia1,weiguo,modulus,lift,cdc15} have studied the discrete-time counterpart.
This paper will focus on the latter
and present a more comprehensive treatment of the work in \cite{cdc15}.
Specifically, the paper provides proofs for the main theorems,
addresses the convergence rate issue, and discusses a time-invariant case
with a less restrictive connectivity condition, which were not included in \cite{cdc15}.

The most general result in the literature regarding the discrete-time version of Altafini's model is that
for any ``repeated jointly strongly connected'' sequence of signed digraphs,
the absolute values of all the agents' opinions will asymptotically reach a consensus,
which is called ``modulus consensus'', having standard consensus and ``bipartite consensus'' 
as special cases \cite{modulus}.
The result is independent of the structure of signs in the digraphs
which can be described by the term structural balance (or structural unbalance) from social sciences \cite{harary}
in that different types of structural balance correspond to different clusterings of opinions in the network
(Section \ref{balance}).
Notwithstanding this, the following questions remain.
What are necessary and sufficient conditions on the sequence of signed digraphs that will lead to a specific clustering?
When will the convergence be exponentially fast and how can the rate of convergence be characterized?
What will happen if the assumption of strong connectivity is relaxed?
This paper aims to answer these questions
and will appeal to a graphical approach prompted by an idea introduced in \cite{lift}, as further discussed below.

In the recent work by Hendrickx \cite{lift},
a very nice lifting approach was proposed to establish the equivalence between
the Altafini model and an expanded DeGroot (or consensus) model with
a special structure; with this equivalence relationship, the convergence results of
the Altafini model were extended under the so-called ``type-symmetry'' assumption
(i.e., the interaction and signs between each pair of neighboring agents are symmetric).
In the discussion section of \cite{lift},
it was mentioned that this lifting approach might ``prove harder to treat systems where
the presence of interaction is symmetric but their signs are not''.
This is precisely what we consider in this paper  and with more generality, we make use of the lifting approach and consider the most general case
where both interaction and signs between each pair of neighboring agents can be asymmetric.
The same approach was used in \cite{weiguo} which provides sufficient conditions
for asymptotic zero consensus and bipartite consensus for a certain class of time-varying
``rooted graphs'' without the ``type-symmetry'' assumption.

The main contributions of this paper are
first, development of a graphical approach to the analysis of the discrete-time version of Altafini's model;
second, derivation of
{\em necessary and sufficient} conditions for {\em exponential} convergence
with respect to different limit states, under a strong connectivity assumption;
third, providing an upper bound on the convergence rate;
and last, description of the system behavior in the case when the neighbor graph is rooted,
but not strongly connected.

The remainder of the paper is organized as follows.
Some notations and preliminaries are introduced in Section \ref{notion}.
In Section \ref{model}, the discrete-time Altafini model is introduced
and the existing modulus consensus result is reviewed.
In Section \ref{balance}, the concepts of structural balance and clustering are introduced.
The main results of the paper are presented in Section \ref{stable},
whose analysis and proofs are given in Section \ref{analysis}.
The convergence rate issue is addressed in Section \ref{rate}.
In Section \ref{discuss}, the system behavior, when the strong connectivity assumption is relaxed, is discussed for
the time-invariant case.
The paper ends with some concluding remarks in Section \ref{ending}.

\subsection{Preliminaries} \label{notion}

For any positive integer $n$, we use $[n]$ to denote the index set $\{1,2,\ldots,n\}$.
We view vectors as column vectors and write $x'$ to denote the transpose of a vector $x$.
For a vector $x$, we use $x_i$ to denote the $i$th entry of $x$.
For any real number $y$, we use $|y|$ to denote its absolute value.
For any matrix $M\in\R^{n\times n}$, we use $m_{ij}$ to denote its $ij$th entry and
write $|M|$ to denote the matrix in $\R^{n\times n}$
whose $ij$th entry is $|m_{ij}|$.
A nonnegative $n\times n$ matrix is called a stochastic matrix if its row sums are all equal to $1$.

For a digraph $\bbb{G}$, we use $(i,j)$ to denote a directed edge from vertex $i$ to vertex $j$.
We say that $\bbb{G}$ has an undirected edge between vertex $i$ and vertex $j$,
denoted by $[i,j]$, if either $(i,j)$ or $(j,i)$ is a directed edge in $\bbb{G}$.
A directed walk of $\bbb{G}$ is a sequence of vertices $v_0,v_1,\ldots,v_m$ in $\bbb{G}$
such that $(v_{i-1},v_{i})$ is a directed edge in $\bbb{G}$ for all $i\in[m]$.
If the vertices $v_0$ and $v_m$ are the same, then the directed walk is called closed.
If all the vertices are distinct, then the directed walk is called a directed path.
If all the vertices are distinct except that vertices $v_0$ and $v_m$ are the same,
then the directed walk is called a directed cycle.

We write $\scr{G}_{sa}$ to denote the set of all digraphs with
  $n$ vertices, 
    which have self-arcs at all vertices.
The graph of an $n\times n$ matrix  $M$ with nonnegative entries is an $n$ vertex directed
 graph $\gamma(M)$  defined so that $(i,j)$ is an arc from $i$ to $j$  in the graph only when
  the $ji$th entry of
 $M$ is nonzero. Such a graph will be in $\scr{G}_{sa}$ if and only if all diagonal entries of
  $M$ are positive.
For purposes of analysis, we write $\bar\scr{G}_{sa}$ to denote the set of all digraphs with
  $2n$ vertices which have self-arcs at all vertices.

A digraph $\mathbb{G}$ is {\em strongly connected} if there is a directed path between
each pair of its distinct vertices. A digraph $\mathbb{G}$ is {\em rooted} if
it contains a directed spanning tree.
A digraph $\mathbb{G}$ is {\em weakly connected} if there is an undirected path between
each pair of its distinct vertices.
Note that every strongly connected graph is rooted and every rooted graph is weakly connected.
The converse statements, however, are false.

\section{The Discrete-Time Altafini Model} \label{model}

Consider a network of $n>1$ agents labeled $1$ to $n$.\footnote{
The purpose of labeling of the agents is only for convenience. We do not require a global labeling of the agents
in the network. We only assume that each agent can identify and differentiate between its neighbors.}
Each agent $i$ has control over a real-valued scalar $x_i(t)$
which the agent is able to update from time to time.
Each agent $i$ can receive information only from certain other agents called agent $i$'s {\em neighbors}.
Neighbor relationships among the $n$ agents are
described by a signed digraph $\bbb{G}(t)$, called {\em neighbor graph}, on $n$ vertices
with an arc from vertex $i$ to vertex $j$ whenever
agent $i$ is a neighbor of agent $j$ at time $t$.
For simplicity, we always take each agent as a neighbor of itself. Thus,
each $\bbb{G}(t)$ has self-arcs at all $n$ vertices.
Each arc is associated with a sign, either positive ``+'' or negative ``--'', which indicates that
agent $i$ regards agent $j$ as a cooperative neighbor if arc $(j,i)$ is associated with a ``+'' sign,
or a competitive neighbor if $(j,i)$ is associated with a ``--'' sign.
It is natural to assume that each self-arc is associated with a ``+'' sign
since an agent cannot compete with itself.

We are interested in the following discrete-time iterative update rule.
At each time $t$ and for each $i\in[n]$, agent $i$ updates its opinion by setting
\eq{x_i(t+1)=\sum_{j\in\scr{N}_i(t)} a_{ij}(t) x_j(t) \label{update}}
where $\scr{N}_i(t)$ denotes the set of neighbors of agent $i$ at time $t$, and each $a_{ij}(t)$ is a real-valued weight
whose sign is consistent with the sign of the arc $(j,i)$.
The weights $a_{ij}(t)$ are assumed to satisfy the following assumption.

\begin{assumption}
For each $i\in[n]$, there hold $a_{ii}(t)>0$
and
$$\sum_{j=1}^n |a_{ij}(t)| = 1 \label{sum}$$
for all time $t$.
There exists a positive number $\beta>0$ such that
$|a_{ij}(t)|\ge \beta$ when $|a_{ij}(t)|>0$ for all $i,j\in[n]$ and $t$.
\label{assume}\end{assumption}

The update rule \rep{update} is an analog of the continuous-time update rule
in \cite{altafini}.
The $n$ update equations in \rep{update} can be combined into one linear recursion equation
\eq{x(t+1)=A(t)x(t) \label{eq}}
in which $x(t)$ is a vector in $\R^n$ whose $i$th entry is $x_i(t)$ and
$A(t)$ is an $n\times n$ matrix whose $ij$th entry is $a_{ij}(t)$.
From Assumption \ref{assume}, it can be seen that the infinite norm of each $A(t)$ equals $1$
and $|A(t)|$ is a stochastic matrix with positive diagonal entries.
In the case when all the arcs have positive signs, the system becomes the standard
linear consensus process which will reach a consensus exponentially fast if and only if
$\bbb{G}(t)$ is ``repeatedly jointly rooted'' \cite{luc,cdc14}.
Thus, system \rep{eq} can be regarded as a generalized model of standard linear consensus.

For each matrix $A(t)$ which satisfies Assumption \ref{assume},
we define the graph of $A(t)$ to be a signed digraph so that
$(i,j)$ is an arc in the graph whenever $a_{ji}(t)$ is nonzero
and the sign of $(i,j)$ is the same as the sign of $a_{ji}(t)$.
It is straightforward to verify that the graph of $A(t)$ is the same as the neighbor graph $\bbb{G}(t)$.
We will use this fact without any special mention in the sequel.

\subsection{Modulus Consensus} \label{modulus}

We say that system \rep{eq} reaches a {\em modulus consensus} if
the absolute values of all $n$ agents $|x_i(t)|$, $i\in[n]$, converge to the same value
as time $t\rightarrow\infty$.
If, in addition, the limiting value does not equal zero and the agents' limiting values have opposite signs,
system \rep{eq} is said to reach a {\em bipartite consensus}.
It should be clear that consensus (including zero consensus and nonzero consensus)
and bipartite consensus are two special cases of modulus consensus.

To state the existing modulus consensus result, we need the following concepts.

We say that a finite sequence of digraphs $\bbb{G}_1,\bbb{G}_2,\ldots,\bbb{G}_m$ with the same vertex set
is {\em jointly strongly connected} if the union\footnote{
The union of a finite sequence of unsigned digraphs with the same vertex set
is an unsigned digraph with the same vertex set and the arc set which is the union
of the arc sets of all digraphs in the sequence.
}
of the digraphs in this sequence is strongly connected.
We say that an infinite  sequence of digraphs $\mathbb{G}_1, \mathbb{G}_2,\ldots$ with the same vertex set
is {\em repeatedly jointly strongly connected}
if there exist positive integers $p$ and $q$ for which
each finite sequence  $\mathbb{G}_{q+kp},\mathbb{G}_{q+kp+1},\ldots,\mathbb{G}_{q+(k+1)p-1}$,
$k\geq 0$, is jointly strongly connected.
It is worth emphasizing that the above connectivity concepts
are also applicable to signed digraphs, without taking signs into account.
Repeatedly jointly strongly connected graphs are sometimes called ``$B$-connected''
graphs in the consensus literature \cite{nedic3}.

The following result states that system \rep{eq} asymptotically reaches
modulus consensus for all initial conditions under an appropriate connectivity assumption.

\begin{proposition}
{\rm (Theorem 2.1 in \cite{modulus})}
Suppose that all $n$ agents adhere to the update rule \rep{update}
and Assumption \ref{assume} holds.
Suppose that the sequence of neighbor graphs $\bbb{G}(1),\bbb{G}(2),\ldots$ is repeatedly jointly strongly connected.
Then, system \rep{eq}
asymptotically reaches a modulus consensus.
\label{ex}\end{proposition}

This proposition was proved in \cite{modulus}, although the form of the connectivity condition is slightly different.
We provide an alternative, more transparent proof in the appendix.

\section{Main Results} \label{main}

As mentioned earlier, modulus consensus includes as special cases bipartite consensus,
zero consensus, and nonzero consensus.
In the sequel, we will explore necessary and sufficient conditions for each type of limit states.
It should be clear that the conditions will naturally depend on the sign structure of
the sequence of neighbor graphs.
We will appeal to the concept of ``structural balance'' from
social sciences \cite{harary}.
The concept of structurally balanced networks was first introduced to consensus problems in \cite{altafini}.
For recent developments on this topic, see \cite{shi15,shi16,bullobalance}.

\subsection{Structural Balance} \label{balance}

A signed digraph $\bbb{G}$ is called {\em structurally balanced} if the vertices of $\bbb{G}$
can be partitioned into two sets such that each arc connecting two agents in the same set
has a positive sign and each arc connecting two agents in different sets
has a negative sign.
Otherwise, the graph $\bbb{G}$ is called {\em structurally unbalanced}.

For our purposes, we extend the same definition of structural balance (and structural unbalance) to ``signed multidigraphs'',
 where a {\em signed multidigraph} is a signed digraph
in which
for some ordered pairs of two distinct vertices $i$ and $j$, there exist two directed edges from vertex $i$ to vertex $j$,
with positive and negative signs respectively.\footnote{
In graph theory, a ``multidigraph'' is usually allowed to have multiple (can be more than two) directed edges
from a vertex $i$ to another vertex $j$. For our purposes, we restrict our attention on
those digraphs which have at most two directed edges from $i$ to $j$. Thus, the set of such $n$-vertex
(signed) digraphs is a finite set.}
Our reason for doing so will be clear shortly.

In this paper, we regard signed multidigraphs as a subset of signed digraphs.
We call those signed digraphs which are not signed multidigraphs, {\em signed simple digraphs}.
In other words, a signed simple digraph does not have multiple directed edges.
Thus, the sets of signed simple digraphs and signed multidigraphs are a partition of
the set of signed digraphs.

There is an equivalent condition for checking structural balance (and structural unbalance) for signed simple digraphs
in the literature \cite{harary,altafini}, as follows.
Let $\bbb{G}$ be a signed simple digraph.
For each directed (or undirected) cycle $\bbb{C}$ in $\bbb{G}$,
we say that $\bbb{C}$ is negative if it contains an odd number of negative signs,
and positive otherwise.
It has been shown that $\bbb{G}$
is structurally balanced if and only if it does not have
negative undirected cycles \cite{harary,altafini}.
Consequently, $\bbb{G}$
is structurally unbalanced if and only if it has at least one
negative undirected cycle.
Using arguments similar to those as in \cite{harary,altafini}, the above equivalent
conditions also hold for signed multidigraphs, and thus hold for all signed digraphs.

%

A simple example of a structurally unbalanced digraph is the signed digraph in which
there exists one pair of distinct agents $i$ and $j$ such that the arcs
$(i,j)$ and $(j,i)$ have different signs.
Another example is the signed multidigraph in which there exist two arcs
from vertex $i$ to another vertex $j$ and they have different signs.
The two graphs in the above examples both have an undirected cycle, consisting of the vertex sequence $i,j,i$,
which is negative.


In the sequel, we will differentiate between different types of structurally balanced signed digraphs
by introducing the concept of clustering.

Let $\scr{I}$ be a set of vectors in $\R^n$ such that
for each $b\in\scr{I}$, there hold $b_1=1$
and $b_i$ equals either $1$ or $-1$ for all $i\in[n]$ and $i\ne 1$.
The set $\scr{I}$ is a finite set and
$\mathbf{1}\in\scr{I}$ where $\mathbf{1}$ denotes the vector in $\R^n$ whose entries all equal $1$.
Each element $b$ in $\scr{I}$ uniquely defines a {\em clustering} of all
the agents in the network by the signs of the entries of $b$.
Specifically, we use $\scr{V}_b^+$ to denote the set of indices in $[n]$
such that $b_i=1$ for all $i\in\scr{V}_b^+$ and
$\scr{V}_b^-$ to denote the set of indices in $[n]$
such that $b_i=-1$ for all $i\in\scr{V}_b^-$.
It can be seen that $\scr{V}_b^+$ and
$\scr{V}_b^-$ are disjoint and $\scr{V}_b^- \cup \scr{V}_b^+ = [n]$.
Since $b_1=1$, it follows that $\scr{V}_b^+$ is nonempty.
In the case when $\scr{V}_b^-$ is nonempty, the vector $b$
defines a unique {\em biclustering} among the agents in the network.
In the special case in which $\scr{V}_b^-$ is an empty set (i.e., $b=\mathbf{1}$),
all the agents belong to the same cluster.

It is worth noting that each possible nonzero modulus consensus (i.e., all the absolute values
of the opinions of the agents in the network reach a consensus at some nonzero value)
can be uniquely represented by an element $b\in\scr{I}$ in that
the agents, including agent $1$, whose labels in  $\scr{V}_b^+$ have the same sign, and
the agents whose labels in  $\scr{V}_b^-$ have the same sign.
In particular, the vector $\mathbf{1}$ represents the standard consensus,
and every other element in $\scr{I}$ represents a unique bipartite consensus.
Thus, all possible nonzero limit states of modulus consensus can be partitioned into
different types, each corresponding to a unique vector in $\scr{I}$.


Each  element $b$ in $\scr{I}$, as well as the unique associated clustering, also
corresponds to a class of signed digraphs with $n$ vertices.
To be more precise,
each element $b\in\scr{I}$ such that $b\ne \mathbf{1}$ corresponds to a class of structurally balanced digraphs,
in which the arcs connecting two vertices in $\scr{V}_b^+$ (or $\scr{V}_b^-$) have positive signs
and the arcs connecting one vertex in $\scr{V}_b^+$ and one vertex in $\scr{V}_b^-$ have negative signs,
and the element $\mathbf{1}$ corresponds to the class of signed digraphs whose signs are all positive.
We call each of the above classes of signed digraphs a {\em structurally balanced class}
and denote it by $\scr{C}_b$, $b\in\scr{I}$.
The remaining signed digraphs with $n$ vertices are all structurally unbalanced and we
call this class of graphs the {\em structurally unbalanced class}, denoted by $\scr{C}_u$.
In general, a signed digraph $\bbb{G}$ may belong to different classes.
But in the case when $\bbb{G}$ is weakly connected, it belongs to a unique class.



\subsection{Exponential Convergence} \label{stable}

For each type of nonzero modulus consensus 
and zero (modulus) consensus, 
we will provide necessary and sufficient conditions under which the modulus consensus can be
reached exponentially fast.
To state our main results, we need the following concepts.


The union of two signed  digraphs $\bbb{G}_p$ and  $\bbb{G}_q$ with the same vertex set
is the signed digraph with the same vertex set,
and the signed arc set being the union of the signed arcs of the two digraphs.
It is clear that the union can be a signed multidigraph.
Note that union is an associative binary operation; because
of this, the definition extends unambiguously to any
finite sequence of signed digraphs including signed multidigraphs.
Since we are interested in signed digraphs with positive self-arcs at all vertices,
the signed digraph generated by the union operation will also have positive self-arcs at all vertices.

We say that a finite sequence of signed digraphs $\bbb{G}_1,\bbb{G}_2,\ldots,\bbb{G}_m$ with the same vertex set
is {\em jointly structurally balanced} with respect to a clustering $b\in\scr{I}$
(or {\em jointly structurally unbalanced}) if the union
of the graphs in this sequence is structurally balanced with respect to the clustering $b$ (or structurally unbalanced).
We say that an infinite  sequence of signed digraphs $\mathbb{G}_1, \mathbb{G}_2,\ldots$ with the same vertex set
is {\em repeatedly jointly structurally balanced} with respect to a clustering $b\in\scr{I}$
(or {\em repeatedly jointly structurally unbalanced})
if there exist positive integers $p$ and $q$ for which
each finite sequence  $\mathbb{G}_{q+kp},\mathbb{G}_{q+kp+1},\ldots,\mathbb{G}_{q+(k+1)p-1}$,
$k\geq 0$, is jointly structurally balanced with respect to the clustering $b$
(or jointly structurally unbalanced).

Note that if a finite sequence of signed digraphs with the same vertex set
is jointly structurally balanced with respect to a clustering $b\in\scr{I}$,
then each graph in the sequence is structurally balanced with respect to the clustering $b$.
Thus, in any sequence of signed digraphs which is jointly structurally balanced with respect to the clustering $b$,
all those signed digraphs which are not in the structurally balanced class $\scr{C}_b$
can appear in the sequence for only a finite number of times.

It is also worth noting that if a finite sequence of signed digraphs $\Sigma$ is jointly structurally unbalanced,
then any finite sequence of signed digraphs which contains $\Sigma$ as a subsequence must be jointly structurally unbalanced.

\begin{remark}
If an infinite  sequence of signed digraphs $\mathbb{G}_1, \mathbb{G}_2,\ldots$ with the same vertex set
is both repeatedly jointly strongly connected and repeatedly jointly structurally balanced with respect to a clustering $b\in\scr{I}$
(or repeatedly jointly structurally unbalanced), then it can be seen that
there must exist positive integers $p$ and $q$ for which
each finite sequence  $\mathbb{G}_{q+kp},\mathbb{G}_{q+kp+1},\ldots,\mathbb{G}_{q+(k+1)p-1}$,
$k\geq 0$, is jointly strongly connected and jointly structurally balanced with respect to the clustering $b$
(or jointly structurally unbalanced).
We will use this
fact without any special mention in the sequel.
\hfill $\Box$
\end{remark}

The main results of this paper are then as follows. We begin with the cases of nonzero modulus consensus.

\begin{theorem}
Suppose that all $n$ agents adhere to the update rule \rep{update}
and Assumption \ref{assume} holds.
Suppose that the sequence of neighbor graphs $\bbb{G}(1),\bbb{G}(2),\ldots$ is repeatedly jointly strongly connected.
Then, for each $b\in\scr{I}$,
system \rep{eq} reaches the corresponding nonzero modulus consensus exponentially fast
for almost all initial conditions
if, and only if, the graph sequence $\bbb{G}(1),\bbb{G}(2),\ldots$
is repeatedly jointly structurally balanced with respect to the clustering $b$.\footnote{
We are indebted to Lili Wang (Department of Electrical Engineering, Yale University)
for pointing out a flaw in the original version of this statement and suggesting how to fix it.
}
\label{nec1}\end{theorem}

The sufficiency of the repeatedly jointly structurally balanced condition
is more or less well known, although the form of the
condition may vary slightly from publication to publication;
see for example \cite{xia,lift,weiguo}.
The proof of necessity will be given in Section \ref{proof}.

\begin{remark}\label{bbb}
Suppose that the sequence of neighbor graphs $\bbb{G}(1),\bbb{G}(2),\ldots$ is repeatedly jointly strongly connected
and structurally balanced with respect to a clustering $b\in\scr{I}$.
Without loss of generality, assume that each $\bbb{G}(t)$, $t\ge 1$,
is structurally balanced with respect to $b$.
Let $B$ be the $n\times n$ diagonal matrix whose $i$th diagonal entry equals $b_i$ for all $i\in[n]$.
Note that $B^2=I$ and for each $A(t)$, the matrix $BA(t)B$ is a stochastic matrix.
Since the sequence of the (signed) graphs of $A(1),A(2),\ldots$ is repeatedly jointly strongly connected,
so is the sequence of the (unsigned) graphs of $BA(1)B,BA(2)B,\ldots$.
With these facts and the well-known result of standard discrete-time linear consensus processes \cite{blondell,luc,ReBe05,reachingp1},
it is straightforward to verify that
the matrix product $A(t)\cdots A(2)A(1)$ converges to a rank one matrix of the form
$bc'$ exponentially fast, where $c$ is a nonzero vector, and thus $x(t)$ converges to $bc'x(1)$.
It follows that system \rep{eq} will reach the corresponding nonzero modulus consensus
if $c'x(1)$ does not equal zero. Since the set of those vectors $x$ which satisfy the equality
$c'x=0$ is a thin set, the nonzero modulus consensus will be reached exponentially fast for almost all initial conditions.
\hfill $\Box$
\end{remark}

\begin{remark}
From the definition of repeatedly jointly structural balance,
Theorem \ref{nec1} implies that
if the system reaches a nonzero modulus consensus corresponding  to a clustering $b$,
then all those signed digraphs which do not belong to the structurally balanced class $\scr{C}_b$
can appear in the sequence of neighbor graphs for only a finite number of times.
\hfill $\Box$
\end{remark}

The next theorem addresses the case of zero (modulus) consensus.

\begin{theorem}
Suppose that all $n$ agents adhere to the update rule \rep{update}
and Assumption \ref{assume} holds.
Suppose that the sequence of neighbor graphs $\bbb{G}(1),\bbb{G}(2),\ldots$ is repeatedly jointly strongly connected.
Then,
system \rep{eq} converges to zero exponentially fast
for all initial conditions
if, and only if, the graph sequence $\bbb{G}(1),\bbb{G}(2),\ldots$
is repeatedly jointly structurally unbalanced.
\label{nec2}\end{theorem}

The proof of this theorem will be given in Section \ref{proof}.

It is worth emphasizing that a repeatedly jointly structurally unbalanced sequence of
signed digraphs does not require each graph in the sequence to be structurally unbalanced.
Actually,  in extreme cases, all graphs in a repeatedly jointly structurally unbalanced sequence can be  structurally balanced,
but with respect to different clusterings.
Let us consider the following example.
Suppose that for each odd time step $2k-1$, $k\ge 1$,
$$A(2k-1)=\matt{0.5 & 0 & 0.5 \cr -0.5 & 0.5 & 0 \cr 0 & -0.5 & 0.5}$$
and for each even time step $2k$, $k\ge 1$,
$$A(2k)=\matt{0.5 & 0 & -0.5 \cr 0.5 & 0.5 & 0 \cr 0 & -0.5 & 0.5}$$
It can be seen that the graphs of $A(2k-1)$ and $A(2k)$ are both structurally balanced,
but the union of the two graphs is structurally unbalanced (and strongly connected).
Thus, according to Theorem \ref{nec2}, system \rep{eq} with $A(2k-1)$ and $A(2k)$
will converge to zero exponentially fast
for all initial conditions.

\begin{remark}
It is well known that for discrete-time linear consensus processes, repeatedly jointly strong connectivity
guarantees exponentially fast consensus \cite{blondell,luc,ReBe05,reachingp1}.
But this is not the case for modulus consensus.
Examples can be generated to show that modulus consensus may be reached only asymptotically, but not exponentially fast.
\hfill $\Box$
\end{remark}

\subsection{Convergence Rate} \label{rate}


Theorems \ref{nec1} and \ref{nec2} provide necessary and sufficient conditions for
system \rep{eq} to converge exponentially fast to different types of modulus consensus.
Of particular interest is the rate at which a modulus consensus is reached.
To characterize the convergence rate, we need the following concept and results.

Let $\{S(t)\}$ be a sequence of stochastic matrices. A sequence of stochastic vectors
$\{\pi(t)\}$ is an absolute probability sequence for $\{S(t)\}$ if
$$\pi'(t) = \pi'(t+1)S(t),\;\;\;\;\; t\geq 1 \label{absolute}$$
It has been shown by Blackwell \cite{blackwell} that every sequence of stochastic matrices
has an absolute probability sequence.
More can be said.

\begin{lemma}
{\rm (Theorem 4.8 in \cite{touri})}
Let  $\{S(t)\}$ be a sequence of $n\times n$ stochastic matrices
that satisfy the following conditions:
\begin{itemize}
\item[(a)] (Aperiodicity) The diagonal entries of each $S(t)$ are positive, i.e., $s_{ii}(t)>0$ for all $t$ and $i\in[n]$.
\item[(b)] (Uniform Positivity) There is a scalar $\beta>0$ such that $s_{ij}(t)\ge\beta$ whenever $s_{ij}(t)>0$.
\item[(c)] (Irreducibility) The sequence of graphs $\{\gamma(S(t))\}$ is repeatedly jointly strongly connected.
\end{itemize}
Let $\{\pi(t)\}$ be  an absolute probability sequence for $\{S(t)\}$.
Then, there is a positive scalar $\delta>0$ such that
$\pi_i(t)\geq \delta$ for all $i\in[n]$ and $t$.
\label{lowbound}
\end{lemma}

A class of stochastic matrices which have this property
was introduced in \cite{touri} (as class $\scr{P}^*$).

The convergence rate of discrete-time linear consensus with a repeatedly jointly strongly connected
sequence of neighbor graphs has been recently characterized with an explicit dependence on the graph
structure including the longest shortest directed path \cite{lya0,lya}.
In the sequel, we will apply the convergence rate result for the discrete-time Altafini model.
We begin with nonzero modulus consensus.

\begin{proposition}
Suppose that all $n$ agents adhere to the update rule \rep{update}
and that Assumption \ref{assume} holds.
Suppose that each neighbor graph $\bbb{G}(t)$ is strongly connected and structurally balanced
with respect to a clustering $b\in\scr{I}$. Then, system \rep{eq} reaches the corresponding
nonzero modulus consensus as fast as $\rho^t$ converges to zero, where
$$\rho = 1- \frac{\delta\beta^2}{4p^*}$$
in which $\delta>0$ is the uniform lower bound on the entries of the absolute probability sequence for
the sequence of stochastic matrices $\{|A(t)|\}$, $\beta>0$ is given in Assumption \ref{assume},
and $p^* = \max_t p^*(t)$ where $p^*(t)$ is the longest shortest directed path of spanning trees contained in $\bbb{G}(t)$.
\label{rate1}\end{proposition}

The proof of this proposition is straightforward using arguments similar to those in
\cite{lya} and is therefore omitted.

We next consider zero (modulus) consensus.

\begin{proposition}
Suppose that all $n$ agents adhere to the update rule \rep{update}
and that Assumption \ref{assume} holds.
Suppose that each neighbor graph $\bbb{G}(t)$ is strongly connected and structurally unbalanced.
Then, system \rep{eq} converges to zero as fast as $\bar\rho^t$ converges to zero, where
$$\bar\rho = 1- \frac{\bar\delta\beta^2}{4\bar p^*}$$
in which $\bar\delta>0$ is the uniform lower bound on the entries of the absolute probability sequence for
the sequence of stochastic matrices $\{\bar A(t)\}$ (matrix $\bar A(t)$ is defined in \rep{large}), $\beta>0$ is given in Assumption \ref{assume},
and $\bar p^* = \max_t \bar p^*(t)$ where $\bar p^*(t)$ is the longest shortest directed path of spanning trees contained in
the graph of $\bar A(t)$.
\label{rate2}\end{proposition}

This proposition is a consequence of Proposition \ref{strong} in Section \ref{graph}.

The results of Propositions \ref{rate1} and \ref{rate2} can be easily extended to
the cases when the sequence of neighbor graphs is repeatedly jointly strongly connected
and structurally balanced (or unbalanced).

\begin{remark}
We will see in the next section that each $n\times n$ matrix $A(t)$ satisfying Assumption \ref{assume}
uniquely determines a $2n\times 2n$ stochastic matrix $\bar A(t)$ defined in \rep{large}.
Thus, $\bar\delta$ and $\bar p^*$ in Proposition \ref{rate2} are the counterparts of
$\delta$ and $p^*$ in Proposition \ref{rate1}.
Using the arguments in the proof of Proposition \ref{strong}, it can be shown that
$\bar p^* \le 2p^* + c^*$ in which
$c^* = \max_t c^*(t)$ where $c^*(t)$ is the longest directed cycle contained in $\bbb{G}(t)$.
Characterization of the relationship between $\delta$ and $\bar\delta$ is a subject of future research.
\hfill $\Box$
\end{remark}

\section{Analysis} \label{analysis}

In this section, we present a graphical approach to analyze the discrete-time
Altafini model.
As mentioned earlier,
the approach we adopt is inspired by a novel idea from \cite{lift}
which lifts the system to an expanded system, as follows.

Define a time-dependent $2n$-dimensional vector $z(t)$ such that for each time $t$,
$$z(t) = \matt{x(t)\cr -x(t)}$$
Then, for all $i\in[2n]$,
$$
z_i(t+1) = \sum_{j=1}^{2n} \bar a_{ij}(t) z_j(t)
$$
in which
\begin{eqnarray*}
\bar a_{ij}(t) \ =\ \bar a_{i+n,j+n}(t) &=& \max\{0,a_{ij}(t)\}  \label{xxx1}\\
\bar a_{i+n,j}(t) \ =\ \bar a_{i,j+n}(t) &=& \max\{0,-a_{ij}(t)\} \label{xxx2}
\end{eqnarray*}
It has been shown in \cite{lift} that the above expanded system is
equivalent to the discrete-time Altafini model.

It is straightforward to verify that the expanded system is a
discrete-time linear consensus process in which the states are coupled.
Thus, it can be written in the form of a state equation
\eq{z(t+1) = \bar A(t)z(t)\label{large}}
where each $\bar A(t)=[\bar a_{ij}(t)]$ is a $2n\times 2n$ stochastic matrix.
With this fact, the graph of $\bar A(t)$ is an unsigned digraph with $2n$ vertices.

The graph of $\bar A(t)$ has the following properties, whose proofs are straightforward and are therefore omitted.

\begin{lemma}
For all $i,j\in[n]$, if $a_{ij}(t)>0$, then the graph of $\bar A(t)$ has an arc from
vertex $j$ to vertex $i$ and an arc from vertex $j+n$ to vertex $i+n$;
if $a_{ij}(t)<0$, then the graph of $\bar A(t)$ has an arc from
vertex $j$ to vertex $i+n$ and an arc from vertex $j+n$ to vertex $i$.
In particular, the graph of $\bar A(t)$ has self-arcs at all $2n$ vertices.
\label{basic}\end{lemma}

\begin{lemma}
Suppose that the graph of $A(t)$ has a directed path from vertex $i$ to vertex $j$ with $i,j\in[n]$.
Then, the graph of $\bar A(t)$ has a directed path from vertex $i$ to vertex $j$ or $j+n$.
In particular, if the directed path from $i$ to $j$ in the graph of $A(t)$ is positive,
then the graph of $\bar A(t)$ has a directed path from $i$ to $j$;
if the directed path from $i$ to $j$ in the graph of $A(t)$ is negative,
then the graph of $\bar A(t)$ has a directed path from $i$ to $j+n$.
\label{path1}\end{lemma}

\begin{lemma}
If the graph of $\bar A(t)$ has a directed path from vertex $i$ to vertex $j$ with $i,j\in[n]$,
then it has a directed path from vertex $i+n$ to vertex $j+n$, and vice versa.
If the graph of $\bar A(t)$ has a directed path from vertex $i$ to vertex $j+n$ with $i,j\in[n]$,
then it has a directed path from vertex $i+n$ to vertex $j$, and vice versa.
\label{path2}\end{lemma}


\subsection{Graphical Results} \label{graph}

In this section, we will establish two key relations between
the signed digraph of $A(t)$ and the expanded unsigned digraph of $\bar A(t)$,
which will play important roles in the proofs of the main results.

\begin{proposition}
Suppose that the graph of $A(t)$ is strongly connected and structurally balanced
with respect to a clustering $b\in\scr{I}$.
Then, the graph of $\bar A(t)$ consists of two disjoint strongly connected components of the same size, $n$.
In particular, the first component consists of vertices $i$, $i\in\scr{V}_b^+$, and $j+n$, $j\in\scr{V}_b^-$,
and the other one consists of vertices $i$, $i\in\scr{V}_b^-$, and $j+n$, $j\in\scr{V}_b^+$.
\label{two}\end{proposition}

This proposition has been proved in \cite{weiguo} (see Lemma 1 in \cite{weiguo}), as well as
the conference version of this paper \cite{cdc15} (see Proposition 1 in \cite{cdc15}).

\begin{proposition}
Suppose that the graph of $A(t)$ is strongly connected and structurally unbalanced.
Then, the graph of $\bar A(t)$ is strongly connected.
\label{strong}\end{proposition}

This result was claimed in \cite{weiguo} (see Lemma 2 in \cite{weiguo}), but with an incomplete\footnote{
In the proof of Lemma 2 in \cite{weiguo},
it is stated that ``since $\bbb{G}$ is structurally unbalanced, there is a negative (directed) cycle
in $\bbb{G}$ \cite{harary}.'' But what \cite{harary} actually proved is that
if $\bbb{G}$ is structurally unbalanced, there is a negative undirected cycle (called semicycle in \cite{harary})
in $\bbb{G}$.}
proof.
In the sequel, we will provide a more complete proof.
The proof relies on the following results.

\begin{lemma}
Suppose that a signed digraph $\bbb{G}$ is strongly connected and structurally unbalanced. Then,
there exists a negative directed closed walk in $\bbb{G}$.
\label{xudong}
\end{lemma}

{\bf Proof:} Suppose that, to the contrary, there does not exist a negative directed walk in $\bbb{G}$.
Since $\bbb{G}$ is structurally unbalanced,
there must exist a negative undirected cycle $\bbb{C}$ in $\bbb{G}$. Label the vertices of  $\bbb{C}$ as $0,1,\ldots, m-1$, with $[0,1], [1,2],\ldots, [m-1,0]$ the associated   undirected edges. In the proof,
we adopt the convention that if an integer $i$ is not in the range $0,1,\ldots, m-1$ but referring to a vertex of $\bbb{C}$, then it refers to the vertex ($i$ mod $m$).

By assumption, the subgraph $\bbb{C}$ cannot be a directed cycle. Thus, there exists at least a vertex $i$ of  $\bbb{C}$ such that
\begin{enumerate}
\item either $i$ is a {\it source}, i.e., both $(i,i-1)$ and $(i,i+1)$ are arcs of $\bbb{G}$,
\item or $i$ is a {\it sink}, i.e., both $(i-1,i)$ and $(i+1,i)$ are arcs of $\bbb{G}$.
\end{enumerate}
Let $\scr{S}$ be the collection of all such vertices. Then, it should be clear that the cardinality of $\scr{S}$ is even. Let $i_1,i_2,\ldots, i_{2k}$ be the elements of $\scr{S}$, with
$$
i_1 < i_2 < \cdots < i_{2k}
$$
Without loss of generality, we assume that $i_1$ is a source. Then, $i_3, i_5,\ldots, i_{2k-1}$ are all sources while $i_2,i_4,\ldots, i_{2k}$ are all sinks. For each $j \in [k]$, let
$$
p^+_{j} := (i_{2j - 1}, i_{2j-1} +1 ) \cdots (i_{2j} -1, i_{2j})
$$
be a directed path of $\bbb{G}$  on $\bbb{C}$ from the source $i_{2j-1}$ to the sink $i_{2j}$. Similarly,  let
$$
p^-_{j} := (i_{2j - 1}, i_{2j-1} - 1 ) \cdots (i_{2j -2} +1, i_{2j -2})
$$
be a directed path of $\bbb{G}$  on $\bbb{C}$ from the source $i_{2j-1}$ to the sink $i_{2j-2}$.

Now fix a $j \in [k]$ and choose a directed path $q^+_{j}$ of $\bbb{G}$ from $i_{2j-2}$ to $i_{2j-1}$. Such a directed path exists since $\bbb{G}$ is strongly connected. Let $n(p^-_j)$ and $n(q^+_j)$ be the numbers of negative signs contained in directed paths $p^-_j$ and $q^+_j$, respectively. Then,
$
n(p^-_j) \equiv n(q^+_j)
$ mod $2$ since otherwise, by concatenating the two directed paths $p^-_j$ and $q^+_j$, we have a negative directed closed walk. Note that this argument applies for all $j$. Now consider a directed closed walk by
concatenating directed paths
$$
w:= q^+_1 p^+_1 q^+_3 p^+_3   \cdots q^+_{2k-1} p^+_{2k-1}
$$
Similarly, we let $n(p^+_j)$ be the number of negative signs contained in $p^+_j$. Then, by the previous arguments, we have
\begin{eqnarray*}
\sum^k_{j=1}\left(n(p^-_{2j-1}) +  n(p^+_{2j-1})\right) \ \ \ \ \ \ \ \ \ \ \ \ \\
\equiv \sum^k_{j=1}\left(n(q^+_{2j-1}) + n(p^+_{2j-1})\right)  \ {\rm mod} \ 2
\end{eqnarray*}
On the other hand, the left hand side of the expression is the total number of negative signs in $\bbb{C}$ which is an odd number. Thus, the right hand side of the expression is also an odd number. In other words, there exists a negative directed closed walk.
\hfill$\qed$

Following Lemma \ref{xudong}, we have the next result.

\begin{corollary}
Suppose that a signed digraph $\bbb{G}$ is strongly connected and structurally unbalanced. Then,
there exists a negative directed cycle  in $\bbb{G}$.
\label{xudong1}\end{corollary}

{\bf Proof:} Let $w$ be the negative directed closed walk in $\bbb{G}$.
Such a directed closed walk exists by Lemma \ref{xudong}.
Choose a vertex in $w$ as a start point and we label it as vertex $i_1$. Express $w$ as
$$
w = (i_1,i_2)(i_2,i_3) \cdots (i_n,i_1)
$$
If $w$ is itself a directed cycle, then the statement is true. Suppose not, then we choose the least integer number $j$ such that
$$
i_j = i_{j+k}
$$
for some $k$. In other words,  $i_j$ is the first vertex, other than $i_1$, in $w$ (with respect to the start point) which appears at least twice  in $w$. For this vertex $i_j$, we may choose the integer $k$ to be the least positive number that $i_j = i_{j+k}$ holds. It should be clear that
$$
(i_j, i_{j+1})\cdots(i_{j+k-1}, i_{j+k})
$$
is a directed cycle. If this directed cycle is negative, then the statement is true. Suppose not, we then remove this positive directed cycle out of the directed closed walk. Then, what remains is still a directed closed walk, denoted by
$$
w':=(i_1,i_2)\cdots (i_{j-1},i_j) (i_j,i_{j+k+1}) \cdots (i_n,i_1)
$$
Moreover, this directed closed walk is also negative. For convenience, we call such an operation a {\it cycle reduction} of a directed closed walk. If $w'$ is a directed cycle, then the statement is true. Suppose not, then we can apply the operation of cycle reduction on $w'$. Thus, we get a sequence of negative directed closed walks as
$$
w \to w^{(1)} \to w^{(2)} \to \cdots
$$
This sequence stops
\begin{enumerate}
\item either at certain step, the removed directed cycle is negative.
\item or there is an integer $l$ such that $w^{(l)}$ is itself a negative directed cycle.
\end{enumerate}
Then, in either of the two cases above, we have found a negative directed cycle. This completes the proof.
\hfill$\qed$

We are now in a position to prove Proposition \ref{strong}.

{\bf Proof of Proposition \ref{strong}:}
By Corollary \ref{xudong1}, there exists a directed negative cycle in the graph of $A(t)$.
For any pair of $i,j\in[n]$, since the graph of $A(t)$ is strongly connected, there
must exist a directed path from vertex $i$ to a vertex $k$ on the directed cycle,
and a directed path from vertex $k$ to vertex $j$.
Thus, there is a directed path from $i$ to $j$, through $k$, in the graph of $A(t)$.
By Lemma \ref{path1}, there exists a directed path in the graph of $\bar A(t)$
from $i$ to $j$ or $j+n$, depending on the sign of the directed path from $i$ to $j$ in the graph of $A(t)$.
Now consider the directed walk from $i$ to $j$ in the graph of $A(t)$ consisting of the above directed path
and a complete round of the directed cycle. Since the directed cycle is negative, the directed walk has a different sign from
the above directed path. Thus,  there exist two directed paths in the graph of $\bar A(t)$
from $i$ to $j$ and to $j+n$. By Lemma \ref{path2}, there exist directed paths
in the graph of $\bar A(t)$ from $i+n$ to $j$ and $j+n$. This completes the proof.
\hfill$\qed$

The results of Propositions \ref{two} and \ref{strong} can be
extended to the cases when a finite sequence of  neighbor graphs is jointly
strongly connected and structurally balanced (or unbalanced), as follows.

\begin{corollary}
Suppose that a finite sequence of the graphs of $A(p),A(p+1),\ldots,A(q)$, $q\ge p$,
is jointly strongly connected and structurally balanced
with respect to a clustering $b\in\scr{I}$.
Then, the union of the graphs of $\bar A(p),\bar A(p+1),\ldots,\bar A(q)$ consists of two disjoint strongly connected components of the same size, $n$.
In particular, the first component consists of vertices $i$, $i\in\scr{V}_b^+$, and $j+n$, $j\in\scr{V}_b^-$,
and the other one consists of vertices $i$, $i\in\scr{V}_b^-$, and $j+n$, $j\in\scr{V}_b^+$.
\label{two1}\end{corollary}

\begin{corollary}
Suppose that a finite sequence of the graphs of $A(p),A(p+1),\ldots,A(q)$, $q\ge p$, is jointly strongly connected and structurally unbalanced.
Then, the union of the graphs of $\bar A(p),\bar A(p+1),\ldots,\bar A(q)$ is strongly connected.
\label{strong1}\end{corollary}

The proofs of Corollaries \ref{two1} and \ref{strong1} are fairly straightforward
using the arguments similar to those used in
the proofs for Propositions \ref{two} and \ref{strong}, and are therefore omitted.

\subsection{Proofs of Main Results} \label{proof}

In this section, we provide proofs for the main results stated in Section \ref{stable}.

{\bf Proof of Theorem \ref{nec1} (Necessity):}
Note that if all those signed digraphs that do not belong to the structurally balanced class $\scr{C}_b$
appear in the sequence of neighbor graphs only a finite number of times, then the sequence of neighbor graphs
is repeatedly jointly structurally balanced with respect to the clustering $b$,
and thus, from Remark \ref{bbb}, the system will reach the corresponding nonzero modulus consensus
exponentially fast for almost all initial conditions.
Thus, to prove the necessity, it is enough to show that
all those signed digraphs that do not belong to $\scr{C}_b$ will only appear a finite number of times.

First consider the signed digraphs in the structurally unbalanced class $\scr{C}_u$.
Since the sequence of neighbor graphs $\bbb{G}(1),\bbb{G}(2),\ldots$
is repeatedly jointly strongly connected, by definition, there exist two positive integers
$p$ and $q$ such that each finite sequence
$\bbb{G}(q+kp),\bbb{G}(q+kp+1),\ldots,\bbb{G}(q+(k+1)p-1)$, $k\ge 0$, is jointly strongly connected.
For each $k\ge 0$, let
$$\bbb{H}_k = \bbb{G}(q+kp)\cup\bbb{G}(q+kp+1)\cup\ldots\cup\bbb{G}(q+(k+1)p-1)$$
Then, each signed digraph $\bbb{H}_k$ is strongly connected.
Note that if a signed digraph $\bbb{G}$ is structurally unbalanced, then any
finite sequence of signed digraphs which contains $\bbb{G}$ must be jointly structurally unbalanced.
Thus, if the graphs in $\scr{C}_u$ appear infinitely many times,
then the graphs in the sequence $\bbb{H}_k$ will be structurally unbalanced for infinitely many times.
By Corollary \ref{strong1}, for each $k\ge 0$,
the union of the graphs of $\bar A(q+kp),\bar A(q+kp+1),\ldots,\bar A(q+(k+1)p-1)$
is strongly connected if $\bbb{H}_k$ is structurally unbalanced.
Thus, the expanded system is a standard discrete-time linear consensus process
whose graph is (jointly) strongly connected infinitely many times.
By Theorem 2 in \cite{hendrickx}, the expanded system will asymptotically reach a consensus.
From the definition of $z(t)$, it follows that
$x(t)$ must converge to zero for all initial conditions.
Thus,
the graphs in $\scr{C}_u$ can only appear in the sequence of neighbor graphs a finite number of times.

Next consider the signed digraphs in the structurally balanced classes.
If only one class of structurally balanced graphs appear infinitely many times and this class is not $\scr{C}_b$,
then from Remark \ref{bbb}, for almost all initial conditions,
the system will reach a nonzero modulus consensus not corresponding to the clustering $b$.
Now suppose that more than one class of structurally balanced graphs appear infinitely many times.
By Theorem 2 in \cite{hendrickx}, each type of modulus consensus corresponding to those structurally balanced classes
will be reached. Since different structurally balanced classes correspond to different types of clusterings,
by the coupled structure of the entries of $z(t)$, the state vector  $x(t)$ must converge to zero for all initial conditions.
Therefore, all the structurally balanced graphs that do not belong to $\scr{C}_b$ can
only appear in the sequence of neighbor graphs a finite number times.
This completes the proof.
\hfill$\qed$

{\bf Proof of Theorem \ref{nec2}:}
We first prove the sufficiency.
Suppose that the sequence of neighbor graphs $\bbb{G}(1),\bbb{G}(2),\ldots$
is repeatedly jointly strongly connected and structurally unbalanced.
Then, there exist two positive integers
$p$ and $q$ such that each finite sequence
$\bbb{G}(q+kp),\bbb{G}(q+kp+1),\ldots,\bbb{G}(q+(k+1)p-1)$, $k\ge 0$, is jointly strongly connected
and structurally unbalanced.
By Corollary \ref{strong1},
the union of the graphs of $\bar A(q+kp),\bar A(q+kp+1),\ldots,\bar A(q+(k+1)p-1)$
is strongly connected for each $k\ge 0$.
Thus, the expanded system is a standard discrete-time linear consensus process
whose sequence of graphs is repeatedly jointly strongly connected.
It is well known that in this case, the expanded vector $z(t)$ will reach a consensus
exponentially fast for all initial conditions \cite{luc,blondell,ReBe05,reachingp1}.
By the definition of $z(t)$, it follows that
$x(t)$ must converge to zero exponentially fast for all initial conditions.

Now we prove the necessity.
Since uniform asymptotic stability and exponential stability are equivalent for linear systems,
it is enough to show that uniform asymptotic stability of system \rep{eq}
implies that the sequence of neighbor graphs
is repeatedly jointly structurally unbalanced.
Suppose therefore that system \rep{eq} is uniformly asymptotically stable.

To establish the claim, suppose that, to the contrary, the sequence of neighbor graphs
$\bbb{G}(1),\bbb{G}(2),\ldots$ is not repeatedly jointly structurally unbalanced. Then,
for every pair of positive integers $l$ and $m$,
there is an integer $k_0> m$ such that the graph sequence
$\bbb{G}(k_0),\bbb{G}(k_0+1),\ldots, \bbb{G}(k_0+l-1)$
is jointly structurally balanced.

Let $\Phi(k,j)$ be the state transition matrix of $A(k)$.
Since $x(k+1)=A(k)x(k)$ is uniformly asymptotically stable,
for each real number $e>0$, there exist
integers $k_e >0$ and $K_e>0$ such that $\|\Phi(k+K_e,k)\|<e$ for all $k>k_e$.
Set $e=1$.
It follows from the preceding arguments that
there must exist an integer $k_0 > k_e$ such that
the graph sequence
$\bbb{G}(k_0),\bbb{G}(k_0+1),\ldots, \bbb{G}(k_0+K_e-1)$
is jointly structurally balanced.

Suppose that the union of the  graphs of
$\bbb{G}(k_0),\bbb{G}(k_0+1),\ldots, \bbb{G}(k_0+K_e-1)$
is structurally balanced with respect to $b\in\scr{I}$.
Let $B$ be the $n\times n$ diagonal matrix whose $i$th diagonal entry equals $b_i$ for all $i\in[n]$.
Then, from Remark \ref{bbb},
\begin{eqnarray*}
\Phi(k_0+K_e,k_0) &=& A(k_0+K_e-1)\cdots A(k_0+1)A(k_0) \\
&=& B\left(BA(k_0+K_e-1)B\right)\cdots \\
&&\left(BA(k_0+1)B\right)\left(BA(k_0)B\right)B
\end{eqnarray*}
in which $BA(i)B$ is a stochastic matrix for all $i\in\{k_0,k_0+1,\ldots,k_0+K_e-1\}$.
It follows that $\Phi(k_0+K_e,k_0)$
has an eigenvalue at $1$ and thus
$\|\Phi(k_0+K_e,k_0)\|=1$, which is a contradiction.
Therefore, the sequence of graphs $\bbb{G}(1),\bbb{G}(2),\ldots$ must be repeatedly jointly structurally unbalanced.
\hfill$\qed$


\section{Discussions} \label{discuss}


For a fixed signed digraph $\bbb{G}$, it has been shown in \cite{xia} that
\begin{enumerate}
  \item If $\bbb{G}$ is structurally balanced and rooted (or strongly connected), then all $|x_i(t)|$ converge to the same value exponentially fast.

  \item If $\bbb{G}$ is structurally unbalanced and strongly connected, then all $x_i(t)$ converge to $0$ exponentially fast.
\end{enumerate}
In this section, we discuss the cases when $\bbb{G}$ is neither structurally balanced nor strongly connected.

Since the least restrictive condition for standard consensus in discrete-time linear
consensus processes is rooted connectivity \cite{luc,cdc14},
we consider the case in which $\bbb{G}$ is structurally unbalanced and rooted, but not necessarily strongly connected.
To state our results, we need the following concept.

We say that vertices $i$ and $j$ in a digraph $\bbb{G}$ are {\em mutually reachable}
if each is reachable from the other along a directed path in $\bbb{G}$.
Mutual reachability is an equivalence relationship on the set of vertices of $\bbb{G}$ which
partitions the vertex set into the disjoint union of a finite number of
equivalence classes.

\begin{lemma}
Each rooted graph has a unique mutually reachable class of roots.
\label{root}
\end{lemma}

The proof of this lemma is fairly straightforward and is therefore omitted.

Lemma \ref{root} implies that if a digraph has more than one root, they are mutually reachable from each other.
Given a rooted digraph $\bbb{G}$. Let $\scr{R}$ denote the set of roots of $\bbb{G}$
and $\bbb{G}_{\scr{R}}$ for the subgraph induced by $\scr{R}$.
From Lemma \ref{root}, $\bbb{G}_{\scr{R}}$ is either a single-vertex graph or
a strongly connected graph.

\begin{proposition}
Suppose that a signed digraph $\bbb{G}$ is structurally unbalanced, rooted, but not strongly connected.
If $\bbb{G}_{\scr{R}}$ is structurally unbalanced, then all the eigenvalues of $A$ are less than $1$ in magnitude.
If $\bbb{G}_{\scr{R}}$ is structurally balanced, then $A$ has an eigenvalue at $1$ and
all the remaining eigenvalues are less than $1$ in magnitude.
\label{new}
\end{proposition}

To prove Proposition \ref{new}, we need the following result.

\begin{lemma}
{\rm (Theorem 8.1.18 in \cite{horn1})}
For any $M\in\R^{n\times n}$,
$\rho\left(M\right) \leq \rho\left(|M|\right)$,
where $\rho(M)$ denotes the spectral radius of matrix $M$.
\label{lili}\end{lemma}

{\bf Proof of Proposition \ref{new}:}
In the case when $\bbb{G}$ is strongly connected, every vertex of $\bbb{G}$ is a root.
Thus, the unique mutually reachable class of roots of $\bbb{G}$ contains all the vertices of $\bbb{G}$.
Since $\bbb{G}$ is not strongly connected, the set of roots of $\bbb{G}$
must be a proper subset of the vertex set of $\bbb{G}$.
Let $\scr{R}=\{i_1,i_2,\ldots,i_m\}$, $m<n$, be the set of roots.
Let $\pi$ be any permutation map for which $\pi(i_j)=j$, $j\in[m]$,
and let $P$ be the corresponding permutation matrix. Then,
$$PAP'=\matt{B & 0 \cr C & D}$$
where $B$ corresponds to the set of roots $\scr{R}$.
Note that $|PAP'|$ is a stochastic matrix with positive diagonal entries.
Then, $|B|$ is also stochastic with positive diagonal entries.
From Lemma \ref{root}, the graph of $B$, $\bbb{G}_{\scr{R}}$, is strongly connected.
It has been proved in \cite{xia} that if $\bbb{G}_{\scr{R}}$ is structurally balanced,
$B$ has an eigenvalue at $1$ and all remaining ones are strictly less than $1$ in magnitude,
and if $\bbb{G}_{\scr{R}}$ is structurally unbalanced, all the eigenvalues of $B$
are strictly less than $1$ in magnitude.

Note that the eigenvalues of $A$ consist of the eigenvalues of $B$ and $D$.
Thus, to prove the lemma, it is sufficient to show that all the eigenvalues of $D$ are
strictly less than $1$ in magnitude.
Suppose the contrary, namely that $D$ has an eigenvalue $\lambda$ whose magnitude is $1$.
Since all the diagonal entries of $A$ are positive,
by the Gersgorin circle theorem, $A$ may have eigenvalues at $1$ and
all other possible eigenvalues must be strictly less than $1$ in magnitude,
so does $PAP'$. Thus, $\lambda=1$.
It is straightforward to verify that
$$P|A|P'=\matt{|B| & 0 \cr |C| & |D|}$$
By Lemma \ref{lili}, $|D|$ must have an eigenvalue at $1$. Note that
$|B|$ is a stochastic matrix, it also has an eigenvalue at $1$.
Thus, $P|A|P'$ has at least two eigenvalues at $1$, so does $|A|$.
But this contradicts the fact that $|A|$ has a unique eigenvalue at $1$ since
$|A|$ is a stochastic matrix with positive diagonal entries whose graph is rooted.
Therefore, all the eigenvalues of $D$ are
strictly less than $1$ in magnitude and the proof is complete.
\hfill$\qed$

More general cases in which the sequence of neighbor graphs is ``repeatedly jointly rooted''
is a subject of future research.

\section{Conclusion} \label{ending}

In this paper, the discrete-time version of Altafini's opinion dynamics model has been
studied through a graphical approach.
Necessary and sufficient conditions for exponential convergence of the system
with respect to different limit states have been established under the assumption of repeatedly jointly strong connectivity.
The rate of convergence has been provided and a time-invariant case without
the assumption of strong connectivity has been discussed.
The time-varying case without the strong connectivity assumption,
which was partially studied in \cite{weiguo}, is a direction for future research.

\section{Acknowledgment}

The authors wish to thank Mahmoud El Chamie (University of Texas at Austin),
Julien M. Hendrickx (Universit\'e Catholique de Louvain),
A. Stephen Morse (Yale University),
Angelia Nedi\'c (University of Illinois at Urbana-Champaign),
Lili Wang (Yale University),
and Zhi Xu (Massachusetts Institute of Technology) for useful discussions
which have contributed to this work.

\bibliographystyle{unsrt}
\bibliography{consensus,social}

\section{Appendix}

{\bf Proof of Proposition \ref{ex}:}
Since the sequence of neighbor graphs $\bbb{G}(1),\bbb{G}(2),\ldots$
is repeatedly jointly strongly connected, by definition, there exist two positive integers
$p$ and $q$ such that each finite sequence
$\bbb{G}(q+kp),\bbb{G}(q+kp+1),\ldots,\bbb{G}(q+(k+1)p-1)$, $k\ge 0$, is jointly strongly connected.
For each $k\ge 0$, let
$$\bbb{H}_k = \bbb{G}(q+kp)\cup\bbb{G}(q+kp+1)\cup\ldots\cup\bbb{G}(q+(k+1)p-1)$$
Then, each signed digraph $\bbb{H}_k$ is strongly connected.
By Corollaries \ref{two1} and \ref{strong1},
the union of the graphs of $\bar A(q+kp),\bar A(q+kp+1),\ldots,\bar A(q+(k+1)p-1)$
is either strongly connected or consists of two strongly connected components.

Note that the set of all $\bbb{H}_k$ is a finite set,
since the number of all possible signed digraphs on $n$ vertices is finite.
This finite set can be partitioned into the sets of structurally balanced graphs $\scr{C}_b$, $b\in\scr{I}$,
and the set of structurally unbalanced graphs $\scr{C}_u$.
It should be clear that there is at least one of these sets whose graphs appear in
the $\bbb{H}_k$ sequence infinitely many times.

Suppose first that structurally unbalanced graphs appear in the $\bbb{H}_k$ sequence infinitely many times.
By Corollary \ref{strong1}, for each $k\ge 0$,
the union of the graphs of $\bar A(q+kp),\bar A(q+kp+1),\ldots,\bar A(q+(k+1)p-1)$
is strongly connected if $\bbb{H}_k$ is structurally unbalanced.
Thus, the expanded system is a standard discrete-time linear consensus process
whose graph is (jointly) strongly connected infinitely many times.
By Theorem 2 in \cite{hendrickx}, the expanded system will asymptotically reach a consensus.
From the definition of $z(t)$, it follows that
$x(t)$ must converge to zero for all initial conditions.

Now suppose that structurally unbalanced graphs appear in the $\bbb{H}_k$ sequence only a finite number of times.
In this case, there exist at least one of the sets $\scr{C}_b$, $b\in\scr{I}$, whose graphs appear infinitely many times.
If only one such set exists, say $\scr{C}_{b_0}$, then from Remark \ref{bbb}, system \rep{eq} will
reach the bipartite consensus with respect to $b_0$ exponentially fast.
If at least two such sets exist, say $\scr{C}_{b_1}$ and $\scr{C}_{b_2}$,
then by Theorem 2 in \cite{hendrickx}, both types of bipartite consensus, with respect to $b_1$ and $b_2$,
will be reached asymptotically. Since different structurally balanced classes correspond to different types of clusterings,
by the coupled structure of the entries of $z(t)$, the state vector  $x(t)$ must converge to zero.

Combining the above cases, we reach the conclusion that system \rep{eq} will always asymptotically reach
a modulus consensus.
\hfill $\qed$

\end{document}